\documentclass[12pt]{article}

\usepackage[left=1.5in, right=1in, top=1in, bottom=1in, dvips,pdftex]{geometry}
\usepackage{color, setspace, graphicx, verbatim, mathrsfs, amsmath, amsthm, amssymb, amsfonts,  multirow, moreverb}
\newtheorem{defn}{Definition}[section]

\newtheorem{thm}[defn]{Theorem}

\theoremstyle{remark}

\numberwithin{equation}{section} \numberwithin{figure}{subsection}

\DeclareMathOperator{\airy}{Ai}

\def\iy{\infty}

\def\be{\begin{equation}}
\def\ee{\end{equation}}

\def\ve{\varepsilon}

\newcommand{\bP}{\mathbb{P}}

\begin{document}

\title{\textbf{On Painlev$\acute{\textbf{e}}$ Related Functions Arising in Random Matrix Theory}}

\author{Leonard N.~Choup \\ Department of Mathematical Sciences \\ University Alabama in Huntsville\\
 Huntsville, AL 35899, USA \\
 email:  \texttt{Leonard.Choup@uah.edu}
} \maketitle

\begin{center} \textbf{Abstract} \end{center}

\begin{small}
In deriving large $n$ probability distribution function of the
rightmost eigenvalue from the classical Random Matrix Theory
Ensembles, one is faced with que question of finding large $n$
asymptotic of certain coupled set of functions. This paper presents
some of these functions in a new light.


\end{small}

\section{Introduction}\label{introduction}
In the study of Edgeworth type expansions for the limiting
distribution of the rightmost eigenvalue from Gaussian Random Matrix
Ensembles, we run into finding large $n$ expansions of many key
functions. For the Tracy-Widom distribution derivation, one needs
the large $n$ limits of these functions, and they can all be express
in terms of the couple pair $q$ and $p$, where $q$ is the
Hastings-McLeod solution to Painlev$\acute{e}$ II equation behaving
at infinity as the Airy function. The frequency of these functions
in the study of the largest eigenvalue of Gaussian and Laguerre
Random Matrix Ensemble points to the necessity of a study of these
functions in their own right. We hope this will shed a light into
understanding some derivations related to this aspect of Random
Matrix Theory and related field making use of such functions.

If one try to read through a proof of an expansion relating various
asymptotic functions, it's easy to get lost in translation. But if
the related functions are well known, the reader will probably have
a different experience and therefore a better understanding of the
techniques and tools used for the derivation.

We present in this paper the derivations of those functions arising
in the study of the largest eigenvalue for the Gaussian Ensemble of
Random Matrix Theory in a hope of achieving our goal set above.

Before stating our results, there is a need to define our
functions.\\

For a Gaussian ensemble of $n \times n$ matrices, the probability
density that the eigenvalues lie in an infinitesimal intervals about
the points $x_{1}\> < \> \ldots \> < \> x_{n}$ is given by
\begin{equation}\label{jpdfeig}
\bP_{n,\beta}(x_{1},\cdots,x_{n})\; = \textrm{C}_{n\beta}\,
\textrm{exp}\left(-\frac{\beta}{2}\, \sum_{1}^{n}x_{j}^{2}\right)\,
\prod_{j<k}|x_{j}-x_{k}|^{\beta}.
\end{equation}
Where $\beta=1$ corresponds to the Gaussian Orthogonal Ensemble
(GOE$_n$),$\beta=2$ corresponds to the Gaussian Unitary Ensemble
(GUE$_{n}$),and $\beta=4$ for the Gaussian Symplectic Ensemble
(GSE$_{n}$). \\
 \eqref{jpdfeig} can also be represented as a determinant involving
 the variables $x_{i}^{'}s$. For the simplest case $\beta=2$, we
 have
 \begin{equation}\label{guedist}
\bP_{n,2}(x_{1},\cdots,x_{n})\; = \frac{1}{n!}(\det
[\varphi_{j-1}(x_{i})])^{2}=\frac{1}{n!}\det[K_{n,2}(x_{i},x_{j})]_{i,j=1,\cdots,n}
\end{equation}
with
\begin{equation}\label{guekernel}
K_{n,2}(x,y)=\sum_{k=0}^{n-1}\varphi_{k}(x)\varphi_{k}(y)=
\sqrt{\frac{n}{2}}\frac{\varphi_{n}(x)\varphi_{n-1}(y)-\varphi_{n-1}(x)\varphi_{n}(y)}{x-y}
\end{equation}
and
\begin{equation*}
 \varphi_{k}(x)= {1\over (2^{k}k! \sqrt{\pi})^{1/2}} \, H_k(x)\,  e^{-x^2/2} \quad \textrm{with}
\quad H_{k}(x) \quad \textrm{the Hermite polinomials}
\end{equation*}

 obtained by orthogonalizing the sequence $\{x^{k},\>
k=0,\cdots,n-1\}$ with respect to $e^{-x^{2}}$ over $\mathbb{R}$.
Using this representation, it can be shown that the probability
distribution function of the largest eigenvalue $\lambda_{max}$ is
given by the Fredholm determinant of the operator with kernel
$K_{n,2}$ acting on the set $(t \>\> \infty)$,

\begin{equation}\label{largeigvaldist}
F_{n,2}(t)=\bP(\lambda_{max} \> < \> t)=\det(I-K_{n,2}).
\end{equation}

In finding the Edgeworth type expansion of $F_{n,2}$, one needs
large $n$ expansion of \eqref{guekernel} or what it amount to, the
large $n$ expansion of $\varphi_{n}$. In \cite{Choup2}, we derived
the following expression.\\
Let the rescaling function $\tau$ be defined by,
\begin{equation}\label{scaling}
\tau (x)=\sqrt{2(n+c)}+ 2^{-\frac{1}{2}}n^{-\frac{1}{6}}x,
\end{equation}
then

\begin{displaymath}
\varphi_{n}(\tau(x))=n^{\frac{1}{6}}  \left\{  \airy(X) +
 \frac{(2c-1)}{2} \airy^{\prime}(X) n^{-\frac{1}{3}} +
 \right.
\left[ (10\,c^{2}-10\,c +\frac{3}{2})\, X \airy(X)\right.
 \end{displaymath}
 \begin{equation}\label{phi(n)}
 +\>\> X^2
 \airy^{\prime}(X)\biggr]\frac{n^{-\frac{2}{3}}}{20}
+ O(n^{-1}) \airy(X) \biggr\}
\end{equation}
and
\begin{displaymath}
\varphi_{n-1}(\tau(x))=n^{\frac{1}{6}}  \left\{  \airy(X) +
 \frac{(2c+1)}{2} \airy^{\prime}(X) n^{-\frac{1}{3}} +
 \right.
\left[ (10\,c^{2}+10\,c +\frac{3}{2})\, X \airy(X) \right.
 \end{displaymath}
 \begin{equation}\label{phi(n-1)}
+ \>\> X^2
 \airy^{\prime}(X)\biggr]\frac{n^{-\frac{2}{3}}}{20} +  O(n^{-1}) \airy(X) \biggr\}
\end{equation}
$Ai$ being the Airy function.
These two functions enable us to obtain the following expansion of the GUE kernel.

\begin{eqnarray*}
K_{n,2}(\tau(X),\tau(Y))\,d\tau(X) =
\tau^{\,_{'}}K_{n}(\tau(X),\tau(Y))dX = \biggl\{ K_{\airy}(X,Y)
-c\airy(X)\airy(Y) n^{-\frac{1}{3}} +
\end{eqnarray*}
\begin{eqnarray*}
\frac{1}{20}\biggl[(X + Y)\airy^{\prime}(X)\airy ^{\prime}(Y) -
(X^2+XY+Y^2)\airy(X)\airy(Y) +
\end{eqnarray*}
\begin{equation}\label{hermitekernel}
 \left. \frac{-20c^2 +3 }{2}(\airy ^{\prime}(X)\airy(Y) +
\airy(X)\airy ^{\prime}(Y)) \right]n^{-\frac{2}{3}}  +O(n^{-1})
E(X,Y)\biggr\}dX.
\end{equation}
In deriving the finite but large $n$ probability distribution function of the largest eigenvalue using \eqref{hermitekernel},  and representation \eqref{largeigvaldist} we have to factor out of \eqref{hermitekernel} the constant term (with respect to $n$) to obtain the representation

\begin{equation*}
F_{n,2}(\tau(t))=\det  \biggl( \left(\,I\>-\>K_{\airy}(X,Y)\right)\cdot\biggl\{ \,I \>
+ \left(\,I\>-\>K_{\airy}(X,Y)\right)^{-1}\biggl[ \> c
\airy(X)\airy(Y) n^{-\frac{1}{3}} -
\end{equation*}
\begin{equation*}
\frac{1}{20}\biggl[(X + Y)\airy^{\prime}(X)\airy ^{\prime}(Y) -
(X^2+XY+Y^2)\airy(X)\airy(Y) +
\end{equation*}
\begin{equation}\label{eq1}
\left. \left.\left. \left.\frac{-20c^2 +3 }{2}(\airy
^{\prime}(X)\airy(Y) + \airy(X)\airy ^{\prime}(Y))
\right]n^{-\frac{2}{3}}  +O(n^{-1}) E(X,Y)
\right]\right\}\right).
\end{equation}
This Fredholm determinant is computed over the set $(t,\> \infty)$.
Thus  to complete the determination of $F_{n,2}(\tau(t))$ we need to determine  the action of the integral operator $(I-K_{\airy})$ on $x^{i}\airy(x)$ and $x^{i}\airy(x)$ where $i=0,1,\cdots$. These are the special functions in the GUE case, and they are independent of $n$, they are well known in the literature (see for example \cite{Choup1,Choup2,Choup3,Trac1,Trac2,Trac3,Trac4,Trac5,Trac7,Trac8}).  For these $n$ independent functions, we just redefine them here and then introduce their $n$ dependent counterparts.

\begin{equation}\label{notation}
K_{\airy}(X,Y)\>=\>
\frac{\airy(X)\,\airy^{'}(Y)\>-\>\airy(Y)\,\airy^{'}(X)}{X-Y}=\int_{0}^{\infty}\airy(X+Z)\,\airy(Y+Z)\,dZ.
\end{equation}
\begin{equation}\label{rho}
\rho(X,Y;s)=(I-K_{\airy})^{-1}(X,Y;s),\quad  R(X,Y;s)=\rho(X,Y;s) \cdot K_{\airy}(X,Y)
\end{equation}
 this last product is operator multiplication.

\begin{equation}\label{Q}
Q_{i}(x;s)\, =\, (\,(I-K_{\airy})^{-1}\, ,\, x^{i}\airy),
\end{equation}
\begin{equation}\label{P}
P_{i}(x;s)\, =\, (\,(I-K_{\airy})^{-1}\, ,\, x^{i}\airy^{'}),
\end{equation}
\begin{equation}\label{q}
q_{i}(s)\>=\> Q_{i}(s;s), \>\>\> p_{i}(s)\>=\> P_{i}(s;s), \>\> q_{0}(s):=q(s)\>\>p_{0}(s):=p(s)
\end{equation}
\begin{equation}\label{u}
u_{i}(s)\,=\, (Q_{i},\airy),\quad v_{i}(s)\,=\, (P_{i},\airy),\>\>u_{0}(s):=u(s),\>\>v_{0}(s):=v(s)
\end{equation}
\begin{equation}
\tilde{v}_{i}(s)=(Q_{i},\airy^{'}),\quad
w_{i}(t)\,= \, (P_{i},\airy^{'}),\>\> w_{0}(s):=w(s),\>\> \textrm{and}\>\>\tilde{v}_{0}(s):=\tilde{v}(s).
\end{equation}
 Here $(\,\cdot \,,\cdot \,)$  denotes
the inner product on $L^{2}(s,\infty)$  and $i=0,1,2,\cdots$. These are all  well known functions, this paper is concerned with the $n$ dependent counterparts whose definitions are similar in nature. the changes needed here are on the kernel definition.  The operator kernel is of the same form as \eqref{hermitekernel}
\begin{equation}\label{guelkernel}
K_{n}(x,y)=\frac{\varphi(x) \psi(y)-\psi(x)\varphi(y)}{x-y}
\end{equation}
with
\begin{equation*}
\varphi(x)=\sqrt[4]{\frac{n}{2}}\>\varphi_{n}(x)\>\>\>\textrm{and}\>\>\> \psi(x)=\sqrt[4]{\frac{n}{2}}\varphi_{n-1}(x)
\end{equation*}
Relating this to the previous set of function are the functions $\varphi$ and $\psi$, they are  $\airy$ and $\airy^{'}$.
We have the following functions,
\begin{equation}\label{rhon}
\rho_{n}(x,y;t):=(I-K_{n})^{-1}(x,y;t), \quad
R_{n}(x,y;t):=\int_{t}^{\infty}\rho_{n}(x,z;t)\> K_{n}(z,y;t)\>dz
\end{equation}
these are kernels of integral operators on $(t \>\> \infty)$
\begin{equation}\label{QniPni}
Q_{n,i}(x;t):=\int_{t}^{\infty}\rho_{n}(x,y;t)y^{i}\varphi(y)\>dy  ,\quad  P_{n,i}(x;t):=\int_{t}^{\infty}\rho_{n}(x,y;t)y^{i}\psi(y) \>dy
\end{equation}
or
\begin{equation*}
Q_{n}(x;t):=(\rho_{n}, \varphi)_{(t \>\> \infty)}\quad P_{n}(x;t):=(\rho_{n},\psi)_{(t,\>\>\infty)}.
\end{equation*}
And the other functions are

\begin{equation}\label{qn}
q_{n,i}(t)\>=\> Q_{n,i}(t;t), \>\>\> p_{n,i}(t)\>=\> P_{n,i}(t;t), \>\> q_{n,0}(t):=q_{n}(t)\>\>p_{n,0}(s):=p_{n}(t)
\end{equation}
\begin{equation}\label{un}
u_{n,i}(t)\,=\, (Q_{n,i},\varphi),\quad v_{n,i}(t)\,=\, (P_{n,i},\varphi),\>\>u_{n,0}(t):=u_{n}(t),\>\>v_{n,0}(t):=v_{n}(t)
\end{equation}
\begin{equation}\label{tildev}
\tilde{v}_{n,i}(t)=(Q_{n,i},\psi),\quad
w_{n,i}(t)\,= \, (P_{n,i},\psi),\>\> w_{n,0}(t):=w_{n}(t),\>\> \textrm{and}\>\>\tilde{v}_{n,0}(t):=\tilde{v}_{n}(t).
\end{equation}
 Here $(\,\cdot \,,\cdot \,)$  denotes
the inner product on $L^{2}(t,\infty)$  and $i=0,1,2,\cdots$.  \\
We will like to point out the following ambiguity in these definitions, the $n$-independent functions have a subscript $i$ whereas the $n$ dependent ones have the subscript $n$. We were not able to find a suitable representations of the set of functions depending on the matrix ensemble of $n\times n$ matrices, but the choice of keeping with the original Tracy and widom notation was made in part to help the reader go through the topic without too much confusion. Thus whenever we use the subscript $n$ we will refer to the large size on the underlying matrix ensemble and when $i$ is used it refers to the exponent of the variable $x$ appearing in the definition of that specific function and $i$ takes values from $0,1,2,\cdots$. One exception is when we will use a second subscript to distinguish between the $3$ beta ensembles $\beta=1,2,4$, in this case we will remind the reader of the significance of those values.\\

In deriving the probability distribution function of the largest eigenvalue $F_{n,1}(t)$ for the orthogonal ensemble, and $F_{n,4}(t)$ for the symplectic ensemble, we encounter new sets on functions obeying the same set of relations.\\
If we define $\ve$ to be  the integral operator with kernel $\varepsilon(x,y)=\frac{1}{2}\textrm{sign} (x-y)$ then
\begin{equation}\label{Qepsilon}
Q_{n,\ve}(x;t):=\int_{t}^{\infty}\rho_{n}(x,y;t)\ve(\varphi)(y)\>dy, \quad q_{n,\ve}(t):=Q_{n,\ve}(t;t)
\end{equation}
\begin{equation}
P_{n,\varepsilon}(x;t):=\int_{t}^{\infty}\rho_{n}(x,y;t)\varepsilon(\psi)(y)\> dy, \quad p_{n,\varepsilon}(t):=P_{n,\varepsilon}(t;t).
\end{equation}
In a similar way we define

\begin{equation}\label{unepsilon}
u_{n,\varepsilon}(t):=\int_{t}^{\infty}Q_{n,\varepsilon}(x;t) \> \varphi(x)\>dx, \quad  v_{n,\varepsilon}(t):=\int_{t}^{\infty}P_{n,\varepsilon}(x;t)\>\varphi(x)\>dx.
\end{equation}
\begin{equation}\label{tildevnepsilon}
\tilde{v}_{n, \varepsilon}(t):=\int_{t}^{\infty}Q_{n,\varepsilon}(x;t)\>\psi(x)\> dx, \quad \textrm{and} \quad  w_{n,\varepsilon}(t)=\int_{t}^{\infty}P_{n,\varepsilon}(x;t)\>\psi(x)\>dx.
\end{equation}
And finally we also have for the Gaussian Orthogonal Ensemble
\begin{equation}\label{calRn1}
\mathcal{R}_{n,1}(t):=\int_{-\infty}^{t} R_{n}(x,t;t)\>dx,\quad \mathcal{P}_{n,1}(t):=\int_{-\infty}^{t} P_{n}(x;t)\>dx,\quad \mathcal{Q}_{n,1}(t):=\int_{-\infty}^{t} Q_{n}(x;t)\>dx
\end{equation}
(Note here that the second subscript here refers to the beta being 1 for the orthogonal ensemble and has nothing to do with the previous discussion on $i$ and $n$.)
For the Gaussian Symplectic Ensemble we have,

\begin{equation}\label{calRn4}
\mathcal{R}_{n,4}(t):=\int_{-\infty}^{\infty}\varepsilon(x,t)R_{n}(x,t;t)\>dx,\quad \mathcal{P}_{n,4}(t)=\int_{-\infty}^{\infty}\varepsilon(x-t)\>P_{n}(x;t)\>dx,\quad \textrm{and} \quad
\end{equation}
\begin{equation*}
 \mathcal{Q}_{n,4}(t):=\int_{-\infty}^{\infty}\varepsilon(x-t)\>Q_{n}(x;t)\>dx,
\end{equation*}
and the $4$ refers to beta being $4$ for the Gaussian Symplectic Ensemble.\\

We have the large $n$ expansion of most of these functions from previous work. What is new in this paper are the large $n$ expansion of $Q_{n,i},\>\>P_{n,i}$ this can be used to derive an expansion for $u_{n,i},\>\>v_{n,i},\>\>\tilde{v}_{n,i},\>\>w_{n,i}$. We also have closed formula for $u_{n,\ve},\>\>\tilde{v}_{n,\ve},\>\>q_{n,\ve},\>\> \mathcal{Q}_{n,1},\>\> \mathcal{P}_{n,1},\>\> \mathcal{R}_{n,1},\>\> \mathcal{Q}_{n,4},\>\> \mathcal{P}_{n,4}$, and $\mathcal{R}_{n,4}$.\\

In the second section we will give a brief justification of $Q_{n,i}$ and $P_{n,i}$   follow in the third section with the justification of these last $9$ functions. Again the motivation for the derivation of these functions is due to their appearance in the Edgeworth type expansion of the largest eigenvalue probability distribution function for the Gaussian Orthogonal and Symplectic Ensembles. \\

\section{Epsilon independent functions}

Building on \eqref{phi(n)}, \eqref{phi(n-1)} and \eqref{hermitekernel} we find that
\begin{displaymath}
Q_{n,i}(x):=((I-K_{n,2})^{-1}(x,y;t),y^{i}\varphi(y))=\int_{t}^{\infty}(I-K_{n,2})^{-1}(x,y;t)\>y^{i}\varphi(y)\>dy
\end{displaymath}
\begin{displaymath}
P_{n,i}(x):=((I-K_{n,2})^{-1}(x,y;t),y^{i}\psi(y))=\int_{t}^{\infty}(I-K_{n,2})^{-1}(x,y;t)\>y^{i}\varphi(y)\>dy
\end{displaymath}
therefore we need to find $\rho_{n}(x,y;t)=(I-K_{n,2})^{-1}(x,y;t)$ in order to find an expressions for these two functions.
But
\begin{equation*}
(I-K_{n,2})^{-1}(\tau(X),\tau(Y);\tau(t))= \biggl\{ \,I \> + \left(\,I\>-\>K_{\airy}\right)^{-1}(X,Y;t)\biggl[ \> c \airy(X)\airy(Y) n^{-\frac{1}{3}} -
\end{equation*}
\begin{equation*}
\frac{1}{20}\biggl[(X + Y)\airy^{\prime}(X)\airy ^{\prime}(Y) - (X^2+XY+Y^2)\airy(X)\airy(Y) +
\end{equation*}
\begin{equation}\label{eq2}
\left. \left. \left.\frac{-20c^2 +3 }{2}(\airy^{\prime}(X)\airy(Y) + \airy(X)\airy ^{\prime}(Y))\right]n^{-\frac{2}{3}}  +O(n^{-1}) E(X,Y)
\right]\right\}^{-1}\cdot\left(\,I\>-\>K_{\airy}(X,Y)\right)^{-1}
\end{equation}
\begin{displaymath}
= \left\{I\>+\>cQ(X)\>\airy (Y) n^{-\frac{1}{3}}- \frac{1}{20}\biggl[ (P_{1}(X)+ Y P(X))\airy ^{\prime}(Y) -
  \right.
\end{displaymath}
\begin{displaymath}
\left.(Q_{2}(X)+Y Q_{1}(X)+Y^{2} Q(X))\>\airy(Y) + \frac{-20c^2 +3 }{2}(P(X)\airy(Y) + Q(X)\airy ^{\prime}(Y))\right]n^{-\frac{2}{3}}
\end{displaymath}
\begin{displaymath}
 \left. \left. +O(\frac{1}{n}) E(X,Y)\right]\right\}^{-1}\cdot\left(\,I\>-\>K_{\airy}(X,Y)\right)^{-1}
\end{displaymath}
\begin{displaymath}
= \left\{I\>-\>cQ(X)\>\airy (Y) n^{-\frac{1}{3}}+ \frac{1}{20}\biggl[ (P_{1}(X)+ Y P(X))\airy ^{\prime}(Y) -
  \right.
\end{displaymath}
\begin{displaymath}
(Q_{2}(X)+Y Q_{1}(X)+Y^{2} Q(X))\>\airy(Y) + \frac{-20c^2 +3 }{2}(P(X)\airy(Y) + Q(X)\airy ^{\prime}(Y))
\end{displaymath}
\begin{displaymath}
 \left. \left. +20c^{2}Q(X;s)u(s)\airy(Y)\biggr]n^{-\frac{2}{3}}
 +O(\frac{1}{n}) E(X,Y)\right]\right\}\cdot\left(\,I\>-\>K_{\airy}(X,Y)\right)^{-1}=\rho_{n}(X,Y;t)
 \end{displaymath}
Note that with this representation of $\rho$ all the $Q_{n,i}$ and $P_{n,i}$ will have no term independent of $n$, but only $Q_{n,0}:=Q_{n}$ and $P_{n,0}:=P_{n}$, in \cite{Choup2} we find that

\begin{equation*}
Q_{n}(\tau(X);\tau(s))=n^{\frac{1}{6}}\biggl[ Q(X;s)+
\left[\frac{2c-1}{2}P(X;s)-c Q(X;s)u(s)\right]n^{-\frac{1}{3}}
\end{equation*}
\begin{equation*}
+\left[(10c^{2}-10c+\frac{3}{2})Q_{1}(X;s)+P_{2}(X;s) +
(-30c^{2}+10c+\frac{3}{2})Q(X;s) v(s) \right.
\end{equation*}
\begin{equation*}
+  P_{1}(X;s) v(s) +P(X;s) v_{1}(s)-Q_{2}(X;s) u(s)-Q_{1}(X;s)
u_{1}(s)-Q(X;s) u_{2}(s)
\end{equation*}
\begin{equation}
+ \left.(-10c^{2}+\frac{3}{2})P(X;s) u(s) +20c^{2}Q(X;s) u^{2}(s)
\right]\frac{n^{-\frac{2}{3}}}{20} +O(n^{-1})E_{q}(X;s)\biggr],
\end{equation}
and
\begin{equation*}
P_{n}(\tau(X);\tau(s))=n^{\frac{1}{6}}\biggl[ Q(X;s)+
\left[\frac{2c+1}{2}P(X;s)-c Q(X;s)u(s)\right]n^{-\frac{1}{3}}
\end{equation*}
\begin{equation*}
+\left[(10c^{2}+10c+\frac{3}{2})Q_{1}(X;s)+P_{2}(X;s) +
(-30c^{2}-10c+\frac{3}{2})Q(X;s) v(s) \right.
\end{equation*}
\begin{equation*}
+  P_{1}(X;s) v(s) +P(X;s) v_{1}(s)-Q_{2}(X;s) u(s)-Q_{1}(X;s)
u_{1}(s)-Q(X;s) u_{2}(s)
\end{equation*}
\begin{equation}
+ \left.(-10c^{2}+\frac{3}{2})P(X;s) u(s) +20c^{2}Q(X;s) u^{2}(s)
\right]\frac{n^{-\frac{2}{3}}}{20} +O(n^{-1})E_{p}(X;s)\biggr].
\end{equation}

Using 

\begin{displaymath}
Q_{n,i}(\tau(X),\tau(s))=(\rho_{n}(\tau(X),\tau(Y),\tau(s)),(\tau(Y))^{i}\varphi(\tau(Y)))_{(\tau(s)\>\>\> \infty)}=
\end{displaymath}
\begin{displaymath}
\sum_{k=0}^{i}\frac{i!}{k! \> (i-k)!}\frac{2^{\frac{i}{2}}(n+c)^{\frac{i-k}{2}}}{2^{k}n^{\frac{k}{2}}}(\rho_{n}(\tau(X),\tau(Y);\tau(s)),Y^{k}\varphi(\tau(Y)))_{(\tau(s)\>\>\> \infty)}.
\end{displaymath}

and
\begin{displaymath}
X^{k}\varphi(\tau(X))=n^{\frac{1}{6}}  \left\{ X^{k} \airy(X) +
 \frac{(2c-1)}{2} X^{k}\airy^{\prime}(X) n^{-\frac{1}{3}} +
 \right.
\left[ (10\,c^{2}-10\,c +\frac{3}{2})\, X^{k+1} \airy(X)\right.
 \end{displaymath}
 \begin{equation}\label{phi(n)}
 +\>\> X^{k+2}
 \airy^{\prime}(X)\biggr]\frac{n^{-\frac{2}{3}}}{20}
+ O(n^{-1}) \airy(X) \biggr\}
\end{equation}
we find that
\begin{displaymath}
\rho(X,Y;s)\cdot \varphi(\tau(X))=n^{\frac{1}{6}}  \left\{  Q_{k}(X;s) +
 \frac{(2c-1)}{2} P_{k}(X) n^{-\frac{1}{3}} +
 \right.
\left[ (10\,c^{2}-10\,c +\frac{3}{2})\, Q_{k+1}(X)\right.
 \end{displaymath}
 \begin{equation}\label{phi(n)}
 +\>\> P_{k+2}(X)\biggr]\frac{n^{-\frac{2}{3}}}{20}
+ O(n^{-1}) Q_{k}(X) \biggr\}.
\end{equation}

Combining this with the action of the first factor on the right of \eqref{eq2} gives the following expression for $(\rho_{n}(\tau(X),\tau(Y);\tau(s)),X^{k}\varphi(\tau(X))$
\begin{displaymath}
n^{\frac{1}{6}}\left\{Q_{k}(X;s)+\left[\frac{2c-1}{2}P_{k}(X;s)-cu_{k}(s)Q(X;s)\right]n^{-\frac{1}{3}}+\left[(10c-20c^{2})\tilde{v}_{k}(s)Q(X;s)+\right.\right.
\end{displaymath}
\begin{displaymath}
 (10c^{2}-10c+\frac{3}{2})Q_{k+1}(X;s)+P_{k+2}(X;s) +P_{1}(X;s)v_{k}(s)+P(X)(Y\airy^{'}(Y), Q_{k}(Y))_{(s\>\> \infty)}
 \end{displaymath}
 \begin{displaymath}
  -u_{k}(s)Q_{2}(X;s)
-Q_{1}(X;s)(Y\airy(Y),Q_{k}(Y;s))_{(s\>\> \infty)}-Q(X;s)(Y^{2}\airy(Y),Q_{k}(Y;s))_{(s\>\> \infty)}
\end{displaymath}
\begin{displaymath}
 \left.+\frac{-20c^{2}+3}{2}P(X;s)u_{k}(s)
+\frac{-20c^{2}+3}{2}Q(X;s)v_{k}(s)+20c^{2}Q(X;s)u(s)u_{k}(s)\right]\frac{n^{-\frac{2}{3}}}{20}
\end{displaymath}
\begin{displaymath}
\left.  +O(\frac{1}{n})E(X;s)\right\}.
\end{displaymath}
To simplify the inner product in this last expression, we use the following recurrence relation derived in \cite{Trac4}\\
$Q_{k}(X;s)=X^{k}Q(X;s)-\sum_{i+j=k-1;i,j>0}(v_{j}Q_{i}-u_{j}P_{i})$ to have
\begin{displaymath}
(Q_{k}(X;s),X\airy(X))_{(s\>\>\> \infty)}=\int_{s}^{\infty}\int_{s}^{\infty}X\airy(X)\rho(X,Y;s)Y^{k}\airy(Y)\>dY\> dX
\end{displaymath}
\begin{displaymath}
=(Q_{1}(X;s),X^{k}\airy(X))=(X Q(X;s)+u(s)P(X;s)-v(s)Q(X;s),X^{k}\airy(X))=
\end{displaymath}
\begin{displaymath}
u_{k+1}(s)+u(s)\tilde{v}_{k}(s)-v(s)u_{k}(s)
\end{displaymath}
and
\begin{displaymath}
(Q_{k}(X;s),X^{2}\airy(X))_{(s\>\>\> \infty)}= (Q_{2}(X;s),X^{k}\airy(X))_{(s\>\>\> \infty)}=
\end{displaymath}
\begin{displaymath}
=(X^{2}Q(X;s)-v(s)(X Q(X;s)+u(s)P(X;s)-v(s)Q(X;s))-u(s)(XP(X;s)
\end{displaymath}
\begin{displaymath}
-w(s)Q(X;s)+v(s)P(X;s))-v_{1}(s)Q(X;s)+u_{1}(s)P(X;s),X^{k}\airy(X))
\end{displaymath}
\begin{displaymath}
=u_{k+2}(s)-v(s)u_{k+1}
-v(s)u(s)\tilde{v}_{k}(s)+v(s)^{2}u_{k}(s)-u(s)\tilde{v}_{k+1}(s)
\end{displaymath}
\begin{displaymath}
+u(s)w(s)u_{k}(s)-u(s)v(s)\tilde{v}_{k}(s)-v_{1}(s)u_{k}(s)+u_{1}(s)\tilde{v}_{k}(s).
\end{displaymath}
we also have

\begin{displaymath}
(Q_{k}(X;s),X\airy^{'}(X))_{(s\>\>\> \infty)}=(P_{1}(X;s),X^{k}\airy(X))=
\end{displaymath}
\begin{displaymath}
\tilde{v}_{k+1}(s)+v(s)\tilde{v}_{k}(s)-w(s)u_{k}(s).
\end{displaymath}
We therefore have

\begin{displaymath}
Q_{n,i}(\tau(X),\tau(s))= \sum_{k=0}^{i}\frac{i!}{k! \> (i-k)!}\frac{2^{\frac{i}{2}-k}(n+c)^{\frac{i-k}{2}}}{n^{\frac{k}{2}-\frac{1}{6}}}.\biggl\{Q_{k}(X;s)+
\end{displaymath}
\begin{displaymath}
\left[\frac{2c-1}{2}P_{k}(X;s)-cu_{k}(s)Q(X;s)
\right]n^{-\frac{1}{3}}+\biggl[(10c-20c^{2})\tilde{v}_{k}(s)Q(X;s)+
\end{displaymath}
\begin{displaymath}
 (10c^{2}-10c+\frac{3}{2})Q_{k+1}(X;s)+P_{k+2}(X;s) +P_{1}(X;s)v_{k}(s)+P(X)\biggl(\tilde{v}_{k+1}(s)+v(s)\tilde{v}_{k}(s)-w(s)u_{k}(s) \biggr)
 \end{displaymath}
 \begin{displaymath}
  -u_{k}(s)Q_{2}(X;s)
-Q_{1}(X;s)\biggl(u_{k+1}(s)+u(s)\tilde{v}_{k}(s)-v(s)u_{k}(s)\biggr)
\end{displaymath}
\begin{displaymath}
-Q(X;s)\biggl(u_{k+2}(s)-v(s)u_{k+1}
-v(s)u(s)\tilde{v}_{k}(s)+v(s)^{2}u_{k}(s)-u(s)\tilde{v}_{k+1}(s)
\end{displaymath}
\begin{displaymath}
+u(s)w(s)u_{k}(s)-u(s)v(s)\tilde{v}_{k}(s)-v_{1}(s)u_{k}(s)+u_{1}(s)\tilde{v}_{k}(s)\biggr)
\end{displaymath}
\begin{displaymath}
 \left.+\frac{-20c^{2}+3}{2}P(X;s)u_{k}(s)
+\frac{-20c^{2}+3}{2}Q(X;s)v_{k}(s)+20c^{2}Q(X;s)u(s)u_{k}(s)\right]\frac{n^{-\frac{2}{3}}}{20}
\end{displaymath}
\begin{equation}\label{Qni}
\left.  +O(\frac{1}{n})E(X;s)\right\}.
\end{equation}

In a similar way we have

\begin{displaymath}
P_{n,i}(\tau(X),\tau(s))=(\rho_{n}(\tau(X),\tau(Y),\tau(s)),(\tau(Y))^{i}\psi(\tau(Y)))_{(\tau(s)\>\>\> \infty)}=
\end{displaymath}
\begin{displaymath}
\sum_{k=0}^{i}\frac{i!}{k! \> (i-k)!}\frac{2^{\frac{i}{2}}(n+c)^{\frac{i-k}{2}}}{2^{k}n^{\frac{k}{2}}}(\rho_{n}(\tau(X),\tau(Y);\tau(s)),Y^{k}\psi(\tau(Y)))_{(\tau(s)\>\>\> \infty)}.
\end{displaymath}
And $(\rho_{n}(\tau(X),\tau(Y);\tau(s)),X^{k}\psi(\tau(X))$ is equal to
\begin{displaymath}
n^{\frac{1}{6}}\left\{Q_{k}(X;s)+\left[\frac{2c+1}{2}P_{k}(X;s)-cu_{k}(s)Q(X;s)\right]n^{-\frac{1}{3}}+\left[-(10c+20c^{2})\tilde{v}_{k}(s)Q(X;s)+\right.\right.
\end{displaymath}
\begin{displaymath}
 (10c^{2}+10c+\frac{3}{2})Q_{k+1}(X;s)+P_{k+2}(X;s) +P_{1}(X;s)v_{k}(s)+P(X)(Y\airy^{'}(Y), Q_{k}(Y))_{(s\>\> \infty)}
 \end{displaymath}
 \begin{displaymath}
  -u_{k}(s)Q_{2}(X;s)
-Q_{1}(X;s)(Y\airy(Y),Q_{k}(Y;s))_{(s\>\> \infty)}-Q(X;s)(Y^{2}\airy(Y),Q_{k}(Y;s))_{(s\>\> \infty)}
\end{displaymath}
\begin{displaymath}
 \left.+\frac{-20c^{2}+3}{2}P(X;s)u_{k}(s)
+\frac{-20c^{2}+3}{2}Q(X;s)v_{k}(s)+20c^{2}Q(X;s)u(s)u_{k}(s)\right]\frac{n^{-\frac{2}{3}}}{20}
\end{displaymath}
\begin{displaymath}
\left.  +O(\frac{1}{n})E(X;s)\right\}.
\end{displaymath}

this therefore gives

\begin{displaymath}
P_{n,i}(\tau(X),\tau(s))=
\sum_{k=0}^{i}\frac{i!}{k! \> (i-k)!}\frac{2^{\frac{i}{2}-k}(n+c)^{\frac{i-k}{2}}}{n^{\frac{k}{2}-\frac{1}{6}}}\times
\end{displaymath}
\begin{displaymath}
\left\{Q_{k}(X;s)+\left[\frac{2c+1}{2}P_{k}(X;s)-cu_{k}(s)Q(X;s)\right]n^{-\frac{1}{3}}+\left[-(10c+20c^{2})\tilde{v}_{k}(s)Q(X;s)+\right.\right.
\end{displaymath}
\begin{displaymath}
 (10c^{2}+10c+\frac{3}{2})Q_{k+1}(X;s)+P_{k+2}(X;s) +P_{1}(X;s)v_{k}(s)+
 \end{displaymath}
 \begin{displaymath}
 P(X)\biggl(\tilde{v}_{k+1}(s)+v(s)\tilde{v}_{k}(s)-w(s)u_{k}(s) \biggr)
 \end{displaymath}
 \begin{displaymath}
  -u_{k}(s)Q_{2}(X;s)
-Q_{1}(X;s)\biggl(u_{k+1}(s)+u(s)\tilde{v}_{k}(s)-v(s)u_{k}(s)\biggr)
\end{displaymath}
\begin{displaymath}
-Q(X;s)\biggl(u_{k+2}(s)-v(s)u_{k+1}
-v(s)u(s)\tilde{v}_{k}(s)+v(s)^{2}u_{k}(s)-u(s)\tilde{v}_{k+1}(s)
\end{displaymath}
\begin{displaymath}
+u(s)w(s)u_{k}(s)-u(s)v(s)\tilde{v}_{k}(s)-v_{1}(s)u_{k}(s)+u_{1}(s)\tilde{v}_{k}(s)\biggr)
\end{displaymath}
\begin{displaymath}
 \left.+\frac{-20c^{2}+3}{2}P(X;s)u_{k}(s)
+\frac{-20c^{2}+3}{2}Q(X;s)v_{k}(s)+20c^{2}Q(X;s)u(s)u_{k}(s)\right]\frac{n^{-\frac{2}{3}}}{20}
\end{displaymath}
\begin{equation}\label{Pni}
\left.  +O(\frac{1}{n})E(X;s)\right\}.
\end{equation}
When we set $i$ to zero we recover $Q_{n}(X;s)$ and $P_{n}(X;s)$.

We see immediately that these two series representations of $Q_{n,i}$ and $P_{n,i}$ are not in terms of $n^{-\frac{1}{3}}$ when $i$ is not zero. We can use \eqref{phi(n)}, \eqref{phi(n-1)}, \eqref{Qni} and \eqref{Pni}, to derive an expansion for $u_{n,i},\>\>v_{n,i},\>\>\tilde{v}_{n,i}$ and $w_{n,i}$ from their representations
\begin{displaymath}
u_{n,i}(t)=(Q_{n,i}(x;t), \varphi(x))_{(t \>\>\> \infty)} \quad v_{n,i}(t)=(P_{n,i}(x;t), \varphi(x))_{(t\>\>\> \infty)}
\end{displaymath}
\begin{displaymath}
\tilde{v}_{n,i}(t)=(Q_{n,i}(x;t),\psi(x))_{(t\>\>\> \infty)}\quad \textrm{and} \quad w_{n,i}(t)=(P_{n,i}(x;t), \psi(x))_{(t\>\>\> \infty)}.
\end{displaymath}
We would like to note that \eqref{Qni} and \eqref{Pni} are the new quantities in this section, as additional corollary the derivation of $q_{n,i}(t)$ and $p_{n,i}(t)$.\\

\begin{displaymath}
q_{n,i}(\tau(s))= \sum_{k=0}^{i}\frac{i!}{k! \> (i-k)!}\frac{2^{\frac{i}{2}-k}(n+c)^{\frac{i-k}{2}}}{n^{\frac{k}{2}-\frac{1}{6}}}.\biggl\{q_{k}(s)+
\end{displaymath}
\begin{displaymath}
\left[\frac{2c-1}{2}p_{k}(s)-cu_{k}(s)q(s)
\right]n^{-\frac{1}{3}}+\biggl[(10c-20c^{2})\tilde{v}_{k}(s)q(s)+
\end{displaymath}
\begin{displaymath}
 (10c^{2}-10c+\frac{3}{2})q_{k+1}(s)+p_{k+2}(s) +p_{1}(s)v_{k}(s)+p(s)\biggl(\tilde{v}_{k+1}(s)+v(s)\tilde{v}_{k}(s)-w(s)u_{k}(s) \biggr)
 \end{displaymath}
 \begin{displaymath}
  -u_{k}(s)q_{2}(s)
-q_{1}(s)\biggl(u_{k+1}(s)+u(s)\tilde{v}_{k}(s)-v(s)u_{k}(s)\biggr)
\end{displaymath}
\begin{displaymath}
-q(s)\biggl(u_{k+2}(s)-v(s)u_{k+1}
-v(s)u(s)\tilde{v}_{k}(s)+v(s)^{2}u_{k}(s)-u(s)\tilde{v}_{k+1}(s)
\end{displaymath}
\begin{displaymath}
+u(s)w(s)u_{k}(s)-u(s)v(s)\tilde{v}_{k}(s)-v_{1}(s)u_{k}(s)+u_{1}(s)\tilde{v}_{k}(s)\biggr)
\end{displaymath}
\begin{displaymath}
 \left.+\frac{-20c^{2}+3}{2}p(s)u_{k}(s)
+\frac{-20c^{2}+3}{2}Q(X;s)v_{k}(s)+20c^{2}q(s)u(s)u_{k}(s)\right]\frac{n^{-\frac{2}{3}}}{20}
\end{displaymath}
\begin{equation}\label{qni}
\left.  +O(\frac{1}{n})e_{q}(s)\right\},
\end{equation}

\begin{displaymath}
p_{n,i}(\tau(s))=
\sum_{k=0}^{i}\frac{i!}{k! \> (i-k)!}\frac{2^{\frac{i}{2}-k}(n+c)^{\frac{i-k}{2}}}{n^{\frac{k}{2}-\frac{1}{6}}}\times
\end{displaymath}
\begin{displaymath}
\left\{q_{k}(s)+\left[\frac{2c+1}{2}p_{k}(s)-cu_{k}(s)q(s)\right]n^{-\frac{1}{3}}+\left[-(10c+20c^{2})\tilde{v}_{k}(s)q(s)+\right.\right.
\end{displaymath}
\begin{displaymath}
 (10c^{2}+10c+\frac{3}{2})q_{k+1}(s)+p_{k+2}(s) +p_{1}(s)v_{k}(s)+
 \end{displaymath}
 \begin{displaymath}
 p(s)\biggl(\tilde{v}_{k+1}(s)+v(s)\tilde{v}_{k}(s)-w(s)u_{k}(s) \biggr)
 \end{displaymath}
 \begin{displaymath}
  -u_{k}(s)q_{2}(s)
-q_{1}(s)\biggl(u_{k+1}(s)+u(s)\tilde{v}_{k}(s)-v(s)u_{k}(s)\biggr)
\end{displaymath}
\begin{displaymath}
-q(s)\biggl(u_{k+2}(s)-v(s)u_{k+1}
-v(s)u(s)\tilde{v}_{k}(s)+v(s)^{2}u_{k}(s)-u(s)\tilde{v}_{k+1}(s)
\end{displaymath}
\begin{displaymath}
+u(s)w(s)u_{k}(s)-u(s)v(s)\tilde{v}_{k}(s)-v_{1}(s)u_{k}(s)+u_{1}(s)\tilde{v}_{k}(s)\biggr)
\end{displaymath}
\begin{displaymath}
 \left.+\frac{-20c^{2}+3}{2}p(s)u_{k}(s)
+\frac{-20c^{2}+3}{2}q(s)v_{k}(s)+20c^{2}q(s)u(s)u_{k}(s)\right]\frac{n^{-\frac{2}{3}}}{20}
\end{displaymath}
\begin{equation}\label{Pni}
\left.  +O(\frac{1}{n})e_{P}(s)\right\}.
\end{equation}

In \cite{Choup2}, we found an expression for  $R_{n}(x,y)=\rho_{n}(x,y) \cdot K_{n,2}(x,y)$, this also follows from \eqref{hermitekernel} and \eqref{eq2}. Note that the following representation will give the same result, \begin{equation}R_{n}(x,y;t)=\frac{Q_{n}(x;t)P_{n}(y;t)-P_{n}(x;t)Q_{n}(y;t)}{x-y}\end{equation}.
\begin{equation*}
R_{n}(\tau(X),\tau(Y);\tau(s))dx= \left[R(X,Y;s)-c\,Q(X;s) Q(Y;s)\,n^{-\frac{1}{3}}\right.
\end{equation*}
\begin{equation*}
+ \frac{n^{-\frac{2}{3}}}{20} \biggl[P_{1}(X;s) P(Y;s) +P(X;s)
P_{1}(Y;s)
\end{equation*}
\begin{equation*}
- Q_{2}(X;s) Q(Y;s) - Q_{1}(X;s) Q_{1}(Y;s) -
 Q(X;s)  Q_{2}(Y;s) + 20 c^{2} u_{0}(s) Q(X;s) Q(Y;s)
\end{equation*}
\begin{equation}\label{eq3}
  +\left. \left. \frac{3-20c^{2}}{2}\left(P(X;s)  Q(Y;s) + Q(X;s)
P(Y;s)\,\right) \right] + O(n^{-1})e_{n}(X,Y)\right]dX.
\end{equation}

\section{Epsilon dependent functions}

The corresponding epsilon functions come from the study of the leftmost eigenvalue from GOE and GSE. We present here the system of equation satisfied by those functions and a solution to these equations leading to our desired functions.\\

To simplify notations we  define
\begin{equation}
V_{n,\varepsilon}(t)=1-\tilde{v}_{n,\varepsilon}(t),\quad
\textrm{and} \quad
\tilde{\mathcal{R}}_{n,1}(t)=1-\mathcal{R}_{n,1}(t).
\end{equation}
With this notation, system  is
\begin{equation}
\frac{d}{dt}\left(
  \begin{array}{c}
    u_{n,\varepsilon}(t) \\
    V_{n,\varepsilon}(t) \\
    q_{n,\varepsilon}(t) \\
  \end{array}
\right) =\left(
               \begin{array}{ccc}
                 0 & 0 & -q_{n}(t) \\
                 0 & 0 & p_{n}(t) \\
                 -p_{n}(t) & q_{n}(t) & 0 \\
               \end{array}
             \right)\,\cdot\,
\left(
  \begin{array}{c}
    u_{n,\varepsilon}(t) \\
    V_{n,\varepsilon}(t) \\
    q_{n,\varepsilon}(t) \\
  \end{array}
\right),
\end{equation}
the boundary conditions in this case are
\begin{equation}
\left(
  \begin{array}{c}
    u_{n,\varepsilon}(\infty) \\
    V_{n,\varepsilon}(\infty) \\
    q_{n,\varepsilon}(\infty) \\
  \end{array}
\right) =\left(
               \begin{array}{c}
               0\\
               1\\
                 c_{\varphi} \\
                 \end{array}
             \right).
\end{equation}

For the orthogonal ensemble
\begin{equation}
\frac{d}{dt}\left(
  \begin{array}{c}
    \mathcal{Q}_{n,1}(t) \\
    \mathcal{P}_{n,1}(t) \\
    \tilde{\mathcal{R}}_{n,1}(t) \\
  \end{array}
\right) =\left(
               \begin{array}{ccc}
                 0 & 0 & q_{n}(t) \\
                 0 & 0 & p_{n}(t) \\
                 p_{n}(t) & q_{n}(t) & 0 \\
               \end{array}
             \right)\,\cdot\,
\left(
  \begin{array}{c}
    \mathcal{Q}_{n,1}(t) \\
    \mathcal{P}_{n,1}(t) \\
    \tilde{\mathcal{R}}_{n,1}(t) \\
  \end{array}
\right),
\end{equation}

with
boundary conditions in this case are
\begin{equation}
\left(
  \begin{array}{c}
    \mathcal{Q}_{n,1}(\infty) \\
    \mathcal{P}_{n,1}(\infty) \\
    \tilde{\mathcal{R}}_{n,1}(\infty) \\
  \end{array}
\right) =\left(
               \begin{array}{c}
                 2c_{\varphi} \\
                 0 \\
                1 \\
               \end{array}
             \right)
\quad \textrm{ as } n \textrm{ is even}.
\end{equation}

We also have for the symplectic ensemble

\begin{equation}
\frac{d}{dt}\left(
  \begin{array}{c}
    \mathcal{Q}_{n,4}(t) \\
    \mathcal{P}_{n,4}(t) \\
    \tilde{\mathcal{R}}_{n,4}(t) \\
  \end{array}
\right) =\left(
               \begin{array}{ccc}
                 0 & 0 & -q_{n}(t) \\
                 0 & 0 & -p_{n}(t) \\
                 -p_{n}(t) & -q_{n}(t) & 0 \\
               \end{array}
             \right)\,\cdot\,
\left(
  \begin{array}{c}
    \mathcal{Q}_{n,4}(t) \\
    \mathcal{P}_{n,4}(t) \\
    \tilde{\mathcal{R}}_{n,4}(t) \\
  \end{array}
\right).
\end{equation}
where $ \tilde{\mathcal{R}}_{n,4}(t) =1+\mathcal{R}_{n,4}(t)$, with corresponding
boundary conditions
\begin{equation}
\left(
  \begin{array}{c}
    \mathcal{Q}_{n,4}(\infty) \\
    \mathcal{P}_{n,4}(\infty) \\
    \tilde{\mathcal{R}}_{n,4}(\infty) \\
  \end{array}
\right) =\left(
               \begin{array}{c}
                 -c_{\varphi} \\
                 -c_{\psi} \\
                1 \\
               \end{array}
             \right)\,=\,
\left(
  \begin{array}{c}
    0 \\
    -c_{\psi} \\
    1 \\
  \end{array}
\right)\quad \textrm{ as } n \textrm{ is odd}.
\end{equation}
The first two set of equations were solved in \cite{Choup3}, here we give the general solution from the series expansion derived there. We will not go back into the derivation, but would like to point out that this is the direct consequence of those matrix exponentials. Our goal here is to give a close formula for those functions.\\
We define
\begin{equation}\label{a(t)}
a(t)=\int_{t}^{\infty}q_{n}(x)\>dx \quad \textrm{and} \quad
b(t)=\int_{t}^{\infty}p_{n}(x)\>dx.
\end{equation}
We note that these two functions scale (under the transformation $\tau$) in the large $n$ limit to the same function \begin{equation*}
\frac{1}{\sqrt{2}}\int_{s}^{\infty}q(x)\> dx\>\>=\>\> \frac{1}{\sqrt{2}}\> \mu(s).
\end{equation*} 
We give this to justify our notation used bellow, and it says that for very large $n$, the argument of all the hyperbolic functions is real. With this notation, we have
\begin{equation}\label{une}
u_{n,\varepsilon}(t)=\frac{a(t)}{2b(t)}[1-\cosh\sqrt{2a(t)b(t)}]+c_{\varphi}\sqrt{\frac{a(t)}{2b(t)}} \sinh\sqrt{2a(t)b(t)},
\end{equation}

\begin{equation}\label{tildevne}
\tilde{V}_{n,\varepsilon}(t)=\frac{1}{2}[1+\cosh\sqrt{2a(t)b(t)}]-c_{\varphi}\sqrt{\frac{b(t)}{2a(t)}} \sinh\sqrt{2a(t)b(t)},
\end{equation}
or
\begin{equation}\label{tildevne}
\tilde{v}_{n,\varepsilon}(t)=1-\frac{1}{2}[1+\cosh\sqrt{2a(t)b(t)}]+c_{\varphi}\sqrt{\frac{b(t)}{2a(t)}} \sinh\sqrt{2a(t)b(t)}
\end{equation}
and
\begin{equation}\label{qne}
q_{n,\varepsilon}(t)=-\sqrt{\frac{a(t)}{2b(t)}} \sinh\sqrt{2a(t)b(t)} +c_{\varphi} \cosh\sqrt{2a(t)b(t)}.
\end{equation}
This result is valid\footnote{The computation for GOE assumes $n$ to be even and for GSE assumes $n$ to be odd} for the GOE when $c_{\varphi} \neq 0$ and for the GSE we have $c_{\varphi}=0$

In the same way we find that for the GOE case, the calligraphic functions are

\begin{equation}\label{calQn1}
\mathcal{Q}_{n,1}(t)=c_{\varphi}[1+\cosh\sqrt{2a(t)b(t)}]-\sqrt{\frac{a(t)}{2b(t)}} \sinh\sqrt{2a(t)b(t)},
\end{equation}

\begin{equation}\label{calPn1}
\mathcal{P}_{n,1}(t)=c_{\varphi}\frac{b(t)}{a(t)}[\cosh\sqrt{2a(t)b(t)}-1]-\sqrt{\frac{b(t)}{2a(t)}} \sinh\sqrt{2a(t)b(t)},
\end{equation}
and
\begin{equation}\label{caltildeRn1}
\tilde{\mathcal{R}}_{n,1}(t)=-2c_{\varphi}\sqrt{\frac{b(t)}{2a(t)}}\sinh\sqrt{2a(t)b(t)}+ \cosh\sqrt{2a(t)b(t)}
\end{equation}
or
\begin{equation}\label{caltildeRn1}
\mathcal{R}_{n,1}(t)=1+2c_{\varphi}\sqrt{\frac{b(t)}{2a(t)}}\sinh\sqrt{2a(t)b(t)}- \cosh\sqrt{2a(t)b(t)}
\end{equation}
where
\begin{equation}\label{c varphi}
c_{\varphi}=(\pi\,n)^{1/4}2^{-3/4 -n/2}\frac{(n!)^{1/2}}{(n/2)!}.
\end{equation}
A large $n$ expansion for $v_{n,\varepsilon}$, $q_{n,\varepsilon}$ is given in \cite{Choup3} on page $17$.  A large $n$ expansion of $\mathcal{P}_{n,1}$ is equation $(3.58)$ and $\mathcal{R}_{n,1}$ is equation $(3.59)$ of the same work. We will therefore give here an expression for $u_{n,\varepsilon}$ and $\mathcal{Q}_{n,1}$ for large $n$.
Substitution of $a(t)=\int_{t}^{\iy}q_{n}(x)\>dx$ and $b(t) =\int_{t}^{\iy}p_{n}(x)\>dx$ into \eqref{une} and \eqref{calQn1} yields the following results.

\begin{thm} \label{large n une}
for $s$ bounded away from minus infinity,
\begin{equation*}
u_{n,\varepsilon}(\tau(s))\>=\> \frac{1}{2}(1-e^{-\mu(s)}) \>+\> \left(\frac{\nu(s)}{4\mu(s)}(e^{-\mu(s)}+\cosh(\mu(s))-2)-\frac{cq(s)}{2}e^{-\mu(s)}\right)n^{-\frac{1}{3}}+
\end{equation*}
\begin{equation*}
\frac{1}{32 \mu(s)^2}\left(e^{-\mu(s)} \left(\nu(s)^2 \left(-\left(-1+e^{\mu(s)}\right) \left(-5-12 c+(3+4
c) e^{\mu(s)}\right)-\right.\right.\right.
\end{equation*}
\begin{equation*}
\left.2 \mu(s) \left(1+6 c-2 c e^{2 \mu(s)}+4 c^2 \mu(s)\right)\right)+
\end{equation*}
\begin{equation*}
4 c \nu(s)\left(3-4 e^{\mu(s)}-e^{2 \mu(s)} (-1+\mu(s))+3 \mu(s)+4 c \mu(s)^2\right) \int_{s}^{\iy} q[x] u[x] \, dx+
\end{equation*}
\begin{equation*}
8 \mu(s) \left(-10 c \left(-3+e^{\mu(s)}\right) \left(-1+e^{\mu(s)}\right) \left(\int_{s}^{\iy} q[x] v[x] \, dx-\int_{s}^{\iy} q_1[x] \,
dx\right)+\right.
\end{equation*}
\begin{equation*}
\mu(s) \left(\left(3-20 c^2\right) \int_{s}^{\iy} p[x] u[x] \, dx+3 \int_{s}^{\iy} q[x] v[x] \, dx+\right.
\end{equation*}
\begin{equation*}
2 \int_{s}^{\iy} v[x] p_1[x] \, dx+2 \int_{s}^{\iy} p_2[x] \, dx+3 \int_{s}^{\iy} q_1[x] \, dx-c^2 \left((\int_{s}^{\iy} q[x] u[x] \, dx)^2-\right.
\end{equation*}
\begin{equation*}
\left.20 \left(2 \int_{s}^{\iy} q[x] u[x]^2 \, dx-3 \int_{s}^{\iy} q[x] v[x] \, dx+\int_{s}^{\iy} q_1[x] \, dx\right)\right)-2 \left(\int_{s}^{\iy} u[x] q_2[x] \, dx+\right.
\end{equation*}
\begin{equation*}
\left.\left.\left.\left.\int_{s}^{\iy} q_1[x] u_1[x] \, dx+\int_{s}^{\iy} q[x] u_2[x] \, dx-\int_{s}^{\iy} p[x] v_1[x] \, dx\right)\right)\right)\right) n^{-\frac{2}{3}} +O(n^{-1})
\end{equation*}
\begin{equation*}
\end{equation*}
We also have

\begin{equation*}
\tilde{v}_{n,\ve}(\tau(s))\>=\>\frac{1}{2} \left(1-e^{-\mu(s)}\right) +(\frac{\nu(s)}{4\mu(s)}\sinh(\mu(s))+\frac{cq(s)}{2}e^{-\mu(s)})n^{-\frac{1}{3}}
\end{equation*}
\begin{equation*}
\frac{1}{16 \mu(s)^2}
\left\{\left(4 c \nu(s) \int_{s}^{\iy} q[x] u[x] \, dx \left(-\cosh[\mu(s)] \mu(s)+2 c e^{-\mu(s)} \mu(s)^2+ \sinh[\mu(s)]\right)\right.\right.
\end{equation*}
\begin{equation*}
+\nu(s)^2 (\cosh[\mu(s)] \mu(s)
\left.\left(-1+4 c-4 c^2 \mu(s)\right)+\left(1-4 c+\mu(s)+4 c^2 \mu(s)^2\right) \sinh[\mu(s)]\right)-
\end{equation*}
\begin{equation*}
4 \mu(s) \left(e^{-\mu(s)} \mu(s) \left(\left(-3+20 c^2\right) \int_{s}^{\iy} p[x] u[x] \, dx-\right.\right.\\
3 \int_{s}^{\iy} q[x] v[x] \, dx-2 \int_{s}^{\iy} v[x] p_1[x] \, dx
\end{equation*}
\begin{equation*}
-2 \int_{s}^{\iy} p_2[x] \, dx-3 \int_{s}^{\iy} q_1[x] \, dx+
c^2 \left((\int_{s}^{\iy} q[x] u[x] \, dx)^2-20 \left(2 \int_{s}^{\iy} q[x] u[x]^2 \, dx \right.\right.
\end{equation*}
\begin{equation*}
\left.\left.-3 \int_{s}^{\iy} q[x] v[x] \, dx+\int_{s}^{\iy} q_1[x] \, dx\right)\right)+
 2 \left(\int_{s}^{\iy} u[x] q_2[x] \, dx+\int_{s}^{\iy} q_1[x] u_1[x] \, dx \right. +
\end{equation*}
\begin{equation*}
\left. \left. \int_{s}^{\iy} q[x] u_2[x] \, dx-\int_{s}^{\iy} p[x] v_1[x] \, dx\right)\right)+\\
\left.\left.20 c \left(\int_{s}^{\iy} q[x] v[x] \, dx-\int_{s}^{\iy} q_1[x] \, dx\right) \sinh[\mu(s)]\right)\right\} n^{-\frac{2}{3}}
\end{equation*}
\begin{equation*}
+O(n^{-1})
\end{equation*}
we also have
\begin{equation*}
q_{n,\varepsilon}(\tau(s))\>=\>\frac{e^{-\mu(s)}}{\sqrt{2}}+(\frac{\nu(s)}{2\sqrt{2}\mu(s)}\sinh\mu(s)+\frac{cq(s)}{\sqrt{2}}e^{-\mu(s)})n^{-\frac{1}{3}}
\end{equation*}
\begin{equation*}
\left\{\frac{1}{8 \sqrt{2} \mu(s)^2}\\
\left(\left(-4 c \nu(s) (\int_{s}^{\iy} q[x] u[x] \, dx) \left(\mu(s) \left(\cosh[\mu(s)]+2 c e^{-\mu(s)} \mu(s)\right)-
\sinh[\mu(s)]\right) \right. \right. \right.+
\end{equation*}
\begin{equation*}
\nu(s)^2 (\cosh[\mu(s)] \mu(s)
\left.\left(1+4 c+4 c^2 \mu(s)\right)-\left(1+4 c+\mu(s)+4 c^2 \mu(s)^2\right) \sinh[\mu(s)]\right)+
\end{equation*}
\begin{equation*}
4 \mu(s) \left(e^{-\mu(s)} \mu(s) \left(\left(-3+20 c^2\right) \int_{s}^{\iy} p[x] u[x] \, dx-\right.\right.\\
3 \int_{s}^{\iy} q[x] v[x] \, dx-2 \int_{s}^{\iy} v[x] p_1[x] \, dx-
\end{equation*}
\begin{equation*}
2 \int_{s}^{\iy} p_2[x] \, dx-3 \int_{s}^{\iy} q_1[x] \, dx+\\
c^2 \left((\int_{s}^{\iy} q[x] u[x] \, dx)^2-20 \left(2 \int_{s}^{\iy} q[x] u[x]^2 \, dx- \right.\right.
\end{equation*}
\begin{equation*}
\left. \left.3 \int_{s}^{\iy} q[x] v[x] \, dx+\int_{s}^{\iy} q_1[x] \, dx\right)\right)+\\
2 \left(\int_{s}^{\iy} u[x] q_2[x] \, dx+\int_{s}^{\iy} q_1[x] u_1[x] \, dx \right.
\end{equation*}
\begin{equation*}
\left. \left.+\int_{s}^{\iy} q[x] u_2[x] \, dx-\int_{s}^{\iy} p[x] v_1[x] \, dx\right)\right)+
\end{equation*}
\begin{equation*}
\left.\left.20 c \left(-\int_{s}^{\iy} q[x] v[x] \, dx+\int_{s}^{\iy} q_1[x] \, dx\right) \sinh[\mu(s)]\right)\right\} n^{-\frac{2}{3}}
\end{equation*}
\end{thm}
and the GOE$_{n}$ calligraphic variables are
\begin{thm} \label{large n Qn1}
for $s$ bounded away from minus infinity,
\begin{equation*}
\mathcal{Q}_{n,1}(\tau(s))\>=\> \frac{1}{\sqrt{2}}(1+e^{-\mu(s)}) \>+\> \left(\frac{\nu(s)}{2\sqrt{2}\mu(s)}\sinh(\mu(s))+\frac{cq(s)}{\sqrt{2}}e^{-\mu(s)}\right)n^{-\frac{1}{3}}+
\end{equation*}
\begin{equation*}
\frac{1}{8 \sqrt{2} \mu(s)^2}
\left(\left(-4 c \nu(s) \int_{s}^{\iy} q[x] u[x] \, dx  \left(\mu(s) \left(\cosh[\mu(s)]+2 c e^{-\mu(s)} \mu(s)\right)-\right.\right.\right.
\end{equation*}
\begin{equation*}
\sinh\mu(s))+\nu(s)^2 (\cosh[\mu(s)] \mu(s) \\
\left.\left(1+4 c+4 c^2 \mu(s)\right)-\left(1+4 c+\mu(s)+4 c^2 \mu(s)^2\right) \sinh\mu(s)\right)+
\end{equation*}
\begin{equation*}
4 \mu(s) \left(e^{-\mu(s)} \mu(s) \left(\left(-3+20 c^2\right) \int_{s}^{\iy} p[x] u[x] \, dx-\right.\right.\\
3 \int_{s}^{\iy} q[x] v[x] \, dx-2 \int_{s}^{\iy} v[x] p_1[x] \, dx
\end{equation*}
\begin{equation*}
-2 \int_{s}^{\iy} p_2[x] \, dx-3 \int_{s}^{\iy} q_1[x] \, dx+
c^2 \left((\int_{s}^{\iy} q[x] u[x] \, dx)^2-20 \left(2 \int_{s}^{\iy} q[x] u[x]^2 \, dx\right. \right.
\end{equation*}
\begin{equation*}
 \left. \left. -3 \int_{s}^{\iy} q[x] v[x] \, dx+\int_{s}^{\iy} q_1[x] \, dx\right)\right)+
\left. 2 \left(\int_{s}^{\iy} u[x] q_2[x] \, dx+\int_{s}^{\iy} q_1[x] u_1[x] \, dx  +\right. \right.
\end{equation*}
\begin{equation*}
 \left. \left. \int_{s}^{\iy} q[x] u_2[x] \, dx-\int_{s}^{\iy} p[x] v_1[x] \, dx\right)\right)+\\
\left.\left.20 c \left(\int_{s}^{\iy} q_1[x] \, dx-\int_{s}^{\iy} q[x] v[x] \, dx\right) \sinh\mu(s)\right)\right) n^{-\frac{2}{3}}
\end{equation*}
\begin{equation*}
+ O(n^{-1}),
\end{equation*}
we also have
\begin{equation*}
\mathcal{P}_{n,1}(\tau(s))\>=\> \left\{\frac{-1+e^{-\mu(s)}}{\sqrt{2}}+ \left(\frac{\nu(s)}{2\sqrt{2} \mu(s)}(e^{-\mu(s)}+\cosh\mu(s))+\frac{cq(s)}{\sqrt{2}}e^{-\mu(s)} \right)n^{-\frac{1}{3}} \right.
\end{equation*}
\begin{equation*}
-
\frac{1}{16 \left(\sqrt{2} \mu(s)^2\right)}\left(\left(e^{-\mu(s)} \left(\nu(s)^2 \left(\left(-1+e^{\mu(s)}\right)
\left(5-12 c+(-3+4 c) e^{\mu(s)}\right)-\right.\right.\right.\right.
\end{equation*}
\begin{equation*}
\left.2 (\mu(s)) \left(1+2 c \left(e^{2 \mu(s)}-3\right)+4 c^2 \mu(s)\right)\right)+4 c \nu(s) \\
\left(4 e^{\mu(s)}-3+e^{2 \mu(s)} (\mu(s)-1)-3 \mu(s)+4 c \mu(s)^2\right)
\end{equation*}
\begin{equation*}
\int_{s}^{\iy} q[x] u[x] \, dx+8 \mu(s) \left(10 c \left(-3+e^{\mu(s)}\right) \left(-1+e^{\mu(s)}\right) \right.\\
\left(\int_{s}^{\iy} q[x] v[x] \, dx-\int_{s}^{\iy} q_1[x] \, dx\right)+
\end{equation*}
\begin{equation*}
\mu(s) \left(\left(3-20 c^2\right) \int_{s}^{\iy} p[x] u[x] \, dx+3 \right.\\
\int_{s}^{\iy} q[x] v[x] \, dx+2 \int_{s}^{\iy} v[x] p_1[x] \, dx+2 \int_{s}^{\iy} p_2[x] \, dx+
\end{equation*}
\begin{equation*}
3 \int_{s}^{\iy} q_1[x] \, dx-c^2 \\
\left((\int_{s}^{\iy} q[x] u[x] \, dx)^2-20 \left(2 \int_{s}^{\iy} q[x] u[x]^2 \, dx-3 \int_{s}^{\iy} q[x] v[x] \, dx + \right.\right.
\end{equation*}
\begin{equation*}
\left. \left.\int_{s}^{\iy} q_1[x] \, dx\right)\right)-
\end{equation*}
\begin{equation*}
\left.\left.\left.\left. 2\left(\int_{s}^{\iy} u[x] q_2[x] \, dx+\int_{s}^{\iy} q_1[x] u_1[x] \, dx+\int_{s}^{\iy} q[x] u_2[x] \, dx-\int_{s}^{\iy} p[x] v_1[x] \, dx\right)\right)\right)\right)\right)
\\
n^{-\frac{2}{3}}
\end{equation*}
\begin{equation*}
+ O(n^{-1})
\end{equation*}
and the last of the GOE$_{n}$ function is\\

\begin{equation*}
\mathcal{R}_{n,1}(\tau(s))\>=\> (1-e^{-\mu(s)}) +\left(\frac{\nu(s)}{2\mu(s)}\sinh\mu(s) -cq(s)e^{-\mu(s)}\right)n^{-\frac{1}{3}} \>\>+
\end{equation*}
\begin{equation*}
\frac{1}{8 \mu(s)^2}\\
\left(\left(4 c \nu(s) \mu(s) \left(-\cosh[\mu(s)] \mu(s)+2 c e^{-\mu(s)} \mu(s)^2+\right.\right.\right.\\
\sinh[\mu(s)])+
\end{equation*}
\begin{equation*}
\nu(s)^2 (\cosh[\mu(s)] \mu(s) \\
\left.\left(-1+4 c-4 c^2 \mu(s)\right)+\left(1-4 c+\mu(s)+4 c^2 \mu(s)^2\right) \sinh[\mu(s)]\right)-
\end{equation*}
\begin{equation*}
4 \mu(s) \left(e^{-\mu(s)} \mu(s) \left(\left(-3+20 c^2\right) \int_{s}^{\iy} p[x] u[x] \, dx-\right.\right.\\
3 \int_{s}^{\iy} q[x] v[x] \, dx-2 \int_{s}^{\iy} v[x] p_1[x] \, dx
\end{equation*}
\begin{equation*}
-2 \int_{s}^{\iy} p_2[x] \, dx-3 \int_{s}^{\iy} q_1[x] \, dx+\\
c^2 \left((\int_{s}^{\iy} q[x] u[x] \, dx)^2-20 \left(2 \int_{s}^{\iy} q[x] u[x]^2 \, dx \right. \right.
\end{equation*}
\begin{equation*}
\left. \left.-3 \int_{s}^{\iy} q[x] v[x] \, dx+\int_{s}^{\iy} q_1[x] \, dx\right)\right)+\\
\left.2 \left(\int_{s}^{\iy} u[x] q_2[x] \, dx+\int_{s}^{\iy} q_1[x] u_1[x] \, dx+ \right. \right.
\end{equation*}
\begin{equation*}
\left. \left. \int_{s}^{\iy} q[x] u_2[x] \, dx-\int_{s}^{\iy} p[x] v_1[x] \, dx\right)\right)+
\end{equation*}
\begin{equation*}
\left.\left.20 c \left(\int_{s}^{\iy} q[x] v[x] \, dx-\int_{s}^{\iy} q_1[x] \, dx\right) \sinh[\mu(s)]\right)\right) n^{-\frac{2}{3}}+ O(n^{-1})
\end{equation*}

\end{thm}

For the symplectic ensemble, the calligraphic variables were mentioned in \cite{Choup3}, the similarity of the corresponding system of differential
equations to the GOE calligraphic system makes our derivation simpler, the coefficient matrices are the same up the minus sign. Keeping with the same notation, we see that
\begin{equation}
\left(
  \begin{array}{c}
    \mathcal{Q}_{n,4}(t) \\
    \mathcal{P}_{n,4}(t) \\
    \tilde{\mathcal{R}}_{n,4}(t) \\
  \end{array}
\right) =\left(
               \begin{array}{ccc}
                 \frac{1}{2}(1+\cosh\sqrt{2ab}) & \frac{a}{2b}(\cosh\sqrt{2ab}-1) & \sqrt{\frac{a}{2b}}\sinh\sqrt{2ab} \\
                 \frac{b}{2a}(\cosh\sqrt{2ab}-1) & \frac{1}{2}(1+\cosh\sqrt{2ab}) & \sqrt{\frac{b}{2a}}\sinh\sqrt{2ab} \\
                 \sqrt{\frac{b}{2a}}\sinh\sqrt{2ab} & \sqrt{\frac{a}{2b}}\sinh\sqrt{2ab}  & \cosh\sqrt{2ab} \\
               \end{array}
             \right)\,\cdot\,
\left(
  \begin{array}{c}
    0 \\
    -c_{\psi} \\
    1 \\
  \end{array}
\right).
\end{equation}
We dropped the $t$ dependence of $a$ and $b$ in the above matrix for esthetic reason.  This gives 

\begin{equation}\label{calQn4}
\mathcal{Q}_{n,4}(t)=-c_{\psi}\frac{a(t)}{2b(t)}[\cosh\sqrt{2a(t)b(t)}-1]-\sqrt{\frac{a(t)}{2b(t)}} \sinh\sqrt{2a(t)b(t)},
\end{equation}

\begin{equation}\label{calPn4}
\mathcal{P}_{n,4}(t)=-c_{\psi}\frac{1}{2}[1+\cosh\sqrt{2a(t)b(t)}]-\sqrt{\frac{b(t)}{2a(t)}} \sinh\sqrt{2a(t)b(t)},
\end{equation}
and
\begin{equation}\label{caltildeRn4}
\tilde{\mathcal{R}}_{n,4}(t)=-c_{\psi}\sqrt{\frac{a(t)}{2b(t)}}\sinh\sqrt{2a(t)b(t)}+ \cosh\sqrt{2a(t)b(t)}
\end{equation}
or
\begin{equation}\label{calRn4}
\mathcal{R}_{n,4}(t)=-c_{\psi}\sqrt{\frac{a(t)}{2b(t)}}\sinh\sqrt{2a(t)b(t)}+ \cosh\sqrt{2a(t)b(t)}-1.
\end{equation}
for the GSE, we have the corresponding formula for $u_{n,\varepsilon}$, $\tilde{v}_{n,\varepsilon}$ and $q_{n,\varepsilon}$.

\begin{thm} \label{large n une gse}
for $s$ bounded away from minus infinity,
\begin{equation*}
u_{n,\varepsilon}(\tau(s))\>=\> -\sinh^{2}(\frac{\mu(s)}{2})+\> \left(\frac{\nu(s)}{2\mu(s)}(\cosh(\mu(s))-1)-\frac{cq(s)}{2}\sinh(\mu(s))\right)n^{-\frac{1}{3}} +
\end{equation*}
\begin{equation*}
\frac{1}{16 \mu(s)^2}\left(\left(8 c \nu(s) \int_{s}^{\iy} q[x] u[x] \, dx \left(-1+\cosh[\mu(s)] \left(1+c \mu(s)^2\right)-\right.\right.\right.
\end{equation*}
\begin{equation*}
\mu(s) \sinh[\mu(s)])+\nu(s)^2 (4+8 c-
\left.4 \cosh[\mu(s)] \left(1+2 c+c^2 \mu(s)^2\right)+(1+8 c) \mu(s) \sinh[\mu(s)]\right)+
\end{equation*}
\begin{equation*}
4 \mu(s) \left(40 c (-1+\cosh[\mu(s)]) \left(-\int_{s}^{\iy} q[x] v[x] \, dx+\int_{s}^{\iy} q_1[x] \, dx\right)+\right.
\end{equation*}
\begin{equation*}
\mu(s) \left(-c^2 \cosh[\mu(s)] (\int_{s}^{\iy} q[x] u[x] \, dx)^2+\left(\left(-3+20 c^2\right) \int_{s}^{\iy} p[x] u[x] \, dx-\right.\right.
\end{equation*}
\begin{equation*}
3 \int_{s}^{\iy} q[x] v[x] \, dx-2 \int_{s}^{\iy} v[x] p_1[x] \, dx-2 \int_{s}^{\iy} p_2[x] \, dx-3 \int_{s}^{\iy} q_1[x] \, dx-
\end{equation*}
\begin{equation*}
20 c^2 \left(2 \int_{s}^{\iy} q[x] u[x]^2 \, dx-3 \int_{s}^{\iy} q[x] v[x] \, dx+\int_{s}^{\iy} q_1[x] \, dx\right)+2 \left(\int_{s}^{\iy} u[x] q_2[x] \, dx+\right.
\end{equation*}
\begin{equation*}
\left.\left.\left.\left.\left.\int_{s}^{\iy} q_1[x] u_1[x] \, dx+\int_{s}^{\iy} q[x] u_2[x] \, dx-\int_{s}^{\iy} p[x] v_1[x] \, dx\right)\right) \sinh[\mu(s)]\right)\right)\right)
\\
n^{-\frac{2}{3}}
\end{equation*}
\begin{equation*}
+O(n^{-1})
\end{equation*}
and
\begin{equation*}
\tilde{v}_{n,\varepsilon}(\tau(s))\>=\> -\sinh^{2}(\frac{\mu(s)}{2})+\> \frac{cq(s)}{2}\sinh(\mu(s))n^{-\frac{1}{3}} +
\end{equation*}
\begin{equation*}
\frac{1}{16 \mu(s)}\left(\left(-4 c^2 \cosh[\mu(s)] \mu(s) (\nu(s)-\int_{s}^{\iy} q[x] u[x] \, dx)^2+\right.\right.
\end{equation*}
\begin{equation*}
\left(\nu(s)^2+4 \mu(s) \left(\left(-3+20 c^2\right) \int_{s}^{\iy} p[x] u[x] \, dx-3 \int_{s}^{\iy} q[x] v[x] \, dx-\right.\right.
\end{equation*}
\begin{equation*}
2 \int_{s}^{\iy} v[x] p_1[x] \, dx-2 \int_{s}^{\iy} p_2[x] \, dx-3 \int_{s}^{\iy} q_1[x] \, dx-20 c^2 \left(2 \int_{s}^{\iy} q[x] u[x]^2 \, dx-\right.
\end{equation*}
\begin{equation*}
\left.3 \int_{s}^{\iy} q[x] v[x] \, dx+\int_{s}^{\iy} q_1[x] \, dx\right)+2 \left(\int_{s}^{\iy} u[x] q_2[x] \, dx+\int_{s}^{\iy} q_1[x] u_1[x] \, dx+\right.
\end{equation*}
\begin{equation*}
\left.\left.\left.\left.\int_{s}^{\iy} q[x] u_2[x] \, dx-\int_{s}^{\iy} p[x] v_1[x] \, dx\right)\right)\right) \sinh[\mu(s)]\right) n^{-\frac{2}{3}}+
O(n^{-1})\\
\end{equation*}
we also have
\begin{equation*}
q_{n,\varepsilon}(\tau(s))\>=\> -\frac{1}{\sqrt{2}}\sinh(\mu(s))+\> \left(\frac{\nu(s)}{2\sqrt{2}\mu(s)}\sinh(\mu(s))+\frac{cq(s)}{\sqrt{2}}\cosh(\mu(s))\right)n^{-\frac{1}{3}} +
\end{equation*}
\begin{equation*}
\frac{1}{8 \sqrt{2} \mu(s)^2}
\left(\left(\cosh[\mu(s)] \mu(s) \left((1+4 c) \nu(s)^2-4 c \nu(s) \int_{s}^{\iy} q[x] u[x] \, dx+\right.\right.\right.
\end{equation*}
\begin{equation*}
4 \mu(s) \left(\left(-3+20 c^2\right) \int_{s}^{\iy} p[x] u[x] \, dx-3 \int_{s}^{\iy} q[x] v[x] \, dx-2 \int_{s}^{\iy} v[x] p_1[x] \, dx-2 \right.
\end{equation*}
\begin{equation*}
\int_{s}^{\iy} p_2[x] \, dx-3 \int_{s}^{\iy} q_1[x] \, dx-20 c^2 \left(2 \int_{s}^{\iy} q[x] u[x]^2 \, dx-3 \int_{s}^{\iy} q[x] v[x] \, dx+\int_{s}^{\iy} q_1[x] \, dx\right)+
\end{equation*}
\begin{equation*}
\left.\left.2 \left(\int_{s}^{\iy} u[x] q_2[x] \, dx+\int_{s}^{\iy} q_1[x] u_1[x] \, dx+\int_{s}^{\iy} q[x] u_2[x] \, dx-\int_{s}^{\iy} p[x] v_1[x] \, dx\right)\right)\right)-
\end{equation*}
\begin{equation*}
\left(\nu(s)^2 \left(1+4 c+4 c^2 \mu(s)^2\right)-4 c \nu(s) \right.\\
\left(1+2 c \mu(s)^2\right) \int_{s}^{\iy} q[x] u[x] \, dx+\\
\end{equation*}
\begin{equation*}
\left.4 c \mu(s) \left(c \mu(s) (\int_{s}^{\iy} q[x] u[x] \, dx)^2+20 \left(\int_{s}^{\iy} q[x] v[x] \, dx-\int_{s}^{\iy} q_1[x] \, dx\right)\right)\right)
\\
\sinh[\mu(s)]) n^{-\frac{2}{3}}
\end{equation*}
\begin{equation*}
+O(n^{-1}).
\end{equation*}
\end{thm}
For the GSE Calligraphic functions we have the following expansions,

\begin{thm} \label{large n une gse}
for $s$ bounded away from minus infinity,
\begin{equation*}
\mathcal{Q}_{n,4}(\tau(s))\>=\frac{1}{2\sqrt{2}}\left(1-\cosh\mu(s)+2\sinh\mu(s)\right)\> +\>
\end{equation*}
\begin{equation*}
\left(\frac{\nu(s)}{2\sqrt{2}\mu(s)}(e^{\mu(s)}-1)-\frac{cq(s)}{2\sqrt{2}}(2\cosh\mu(s)-\sinh\mu(s))\right)n^{-\frac{1}{3}} +
\end{equation*}
\begin{equation*}
\frac{1}{32 \sqrt{2} \mu(s)^2}\left(e^{-\mu(s)} \left(\nu(s)^2 \left(-2 \left(-1+e^{\mu(s)}\right) \left(-3-8
c+e^{\mu(s)}\right)-\right.\right.\right.\\
\end{equation*}
\begin{equation*}
\left.\left(3+16 c+e^{2 \mu(s)}\right) \mu(s)+4 c^2 \left(-3+e^{2 \mu(s)}\right) \mu(s)^2\right)-8 c \nu(s) \\
\end{equation*}
\begin{equation*}
\left(2 \left(-1+e^{\mu(s)}\right)+\mu(s) \left(-2+c \left(-3+e^{2 \mu(s)}\right) \mu(s)\right)\right) \int_{s}^{\iy}
q[x] u[x] \, dx+\\
\end{equation*}
\begin{equation*}
4 \mu(s) \left(80 c \left(-1+e^{\mu(s)}\right) \left(\int_{s}^{\iy} q[x] v[x] \, dx-\int_{s}^{\iy} q_1[x] \, dx\right)+\mu(s) \right.\\
\end{equation*}
\begin{equation*}
\left(-\left(-3+20 c^2\right) \left(3+e^{2 \mu(s)}\right) \int_{s}^{\iy} p[x] u[x] \, dx+c^2 \left(-3+e^{2 \mu(s)}\right) (\int_{s}^{\iy} q[x] u[x]
\, dx)^2+\right.\\
\end{equation*}
\begin{equation*}
\left(3+e^{2 \mu(s)}\right) \left(40 c^2 \int_{s}^{\iy} q[x] u[x]^2 \, dx+\left(3-60 c^2\right) \int_{s}^{\iy} q[x] v[x] \, dx+\right.\\
\end{equation*}
\begin{equation*}
2 \int_{s}^{\iy} v[x] p_1[x] \, dx+2 \int_{s}^{\iy} p_2[x] \, dx+\left(3+20 c^2\right) \int_{s}^{\iy} q_1[x] \, dx-2 \left(\int_{s}^{\iy} u[x] q_2[x] \, dx+\right.\\
\end{equation*}
\begin{equation*}
\left.\left.\left.\left.\left.\int_{s}^{\iy} q_1[x] u_1[x] \, dx+\int_{s}^{\iy} q[x] u_2[x] \, dx-\int_{s}^{\iy} p[x] v_1[x] \, dx\right)\right)\right)\right)\right) n^{-\frac{2}{3}} +O(n^{-1})
\end{equation*}
the next function is
\begin{equation*}
\mathcal{P}_{n,4}(\tau(s))\>=\> \frac{1}{2\sqrt{2}}(2\sinh\mu(s)-\cosh\mu(s)-1)
\end{equation*}
\begin{equation*}
+\left(\frac{\nu(s)}{2\sqrt{2}\mu(s)}\sinh\mu(s)-\frac{cq(s)}{2\sqrt{2}}\left(2\cosh\mu(s)-\sinh\mu(s)\right)\right)n^{-\frac{1}{3}}
\end{equation*}
\begin{equation*}
\frac{1}{16 \sqrt{2} \mu(s)^2}((8 c \nu(s) \int_{s}^{\iy} q[x] u[x] \, dx (-\cosh[\mu(s)] \mu(s)+\\
\end{equation*}
\begin{equation*}
\left.c \mu(s)^2 (\cosh[\mu(s)]-2 \sinh[\mu(s)])+\sinh[\mu(s)]\right)+\\
\end{equation*}
\begin{equation*}
\nu(s)^2 \left(-2 \cosh[\mu(s)] \mu(s) \left(1-4 c+2 c^2 \mu(s)\right)+\right.\\
\end{equation*}
\begin{equation*}
\left.\left(2-8 c+\mu(s)+8 c^2 \mu(s)^2\right) \sinh[\mu(s)]\right)-\\
\end{equation*}
\begin{equation*}
4 \mu(s) \left(\mu(s) \left(c^2 (\int_{s}^{\iy} q[x] u[x] \, dx)^2 (\cosh[\mu(s)]-2 \sinh[\mu(s)])+\right.\right.\\
\end{equation*}
\begin{equation*}
\left(-3+20 c^2\right) (\int_{s}^{\iy} p[x] u[x] \, dx) (2 \cosh[\mu(s)]-\sinh[\mu(s)])-\\
\end{equation*}
\begin{equation*}
\left(40 c^2 \int_{s}^{\iy} q[x] u[x]^2 \, dx+\left(3-60 c^2\right) \int_{s}^{\iy} q[x] v[x] \, dx+2 \int_{s}^{\iy} v[x] p_1[x] \, dx+\right.\\
\end{equation*}
\begin{equation*}
2 \int_{s}^{\iy} p_2[x] \, dx+\left(3+20 c^2\right) \int_{s}^{\iy} q_1[x] \, dx-2 \left(\int_{s}^{\iy} u[x] q_2[x] \, dx+\int_{s}^{\iy} q_1[x] u_1[x] \, dx+\right.\\
\end{equation*}
\begin{equation*}
\left.\left.\left.\int_{s}^{\iy} q[x] u_2[x] \, dx-\int_{s}^{\iy} p[x] v_1[x] \, dx\right)\right) (2 \cosh[\mu(s)]-\sinh[\mu(s)])\right)+\\
\end{equation*}
\begin{equation*}
\left.\left.40 c \left(\int_{s}^{\iy} q[x] v[x] \, dx-\int_{s}^{\iy} q_1[x] \, dx\right) \sinh[\mu(s)]\right)\right) n^{-\frac{2}{3}}+O(n^{-1})
\end{equation*}
the last of our function is
\begin{equation*}
\mathcal{R}_{n,4}(\tau(s))\>=\> \left( \cosh\mu(s)-\frac{1}{2}\sinh\mu(s) -1\right)+
\end{equation*}
\begin{equation*}
\left(\frac{\nu(s)}{4\mu(s)}\sinh\mu(s) +\frac{cq(s)}{2}\left(\cosh\mu(s)-2\sinh\mu(s)\right)\right)n^{-\frac{1}{3}}
\end{equation*}
\begin{equation*}
\frac{1}{16 \mu(s)^2}((-4 c \nu(s) (\int_{s}^{\iy} q[x] u[x] \, dx) (\cosh[\mu(s)] \mu(s) \\
\end{equation*}
\begin{equation*}
\left.(1+4 c \mu(s))-\left(1+2 c \mu(s)^2\right) \sinh[\mu(s)]\right)+\\
\end{equation*}
\begin{equation*}
\nu(s)^2 \left(\cosh[\mu(s)] \mu(s) \left(1+4 c+8 c^2 \mu(s)\right)-\right.\\
\end{equation*}
\begin{equation*}
\left.\left(1+4 c+2 \mu(s)+4 c^2 \mu(s)^2\right) \sinh[\mu(s)]\right)+4 \mu(s) \\
\end{equation*}
\begin{equation*}
\left(\mu(s) \left(\left(-3+20 c^2\right) \int_{s}^{\iy} p[x] u[x] \, dx (\cosh[\mu(s)]-2 \sinh[\mu(s)])-\right.\right.\\
\end{equation*}
\begin{equation*}
\left(40 c^2 \int_{s}^{\iy} q[x] u[x]^2 \, dx+\left(3-60 c^2\right) \int_{s}^{\iy} q[x] v[x] \, dx+2 \int_{s}^{\iy} v[x] p_1[x] \, dx+\right.\\
\end{equation*}
\begin{equation*}
2 \int_{s}^{\iy} p_2[x] \, dx+\left(3+20 c^2\right) \int_{s}^{\iy} q_1[x] \, dx-2 \left(\int_{s}^{\iy} u[x] q_2[x] \, dx+\int_{s}^{\iy} q_1[x] u_1[x] \, dx+\right.\\
\end{equation*}
\begin{equation*}
\left.\left.\int_{s}^{\iy} q[x] u_2[x] \, dx-\int_{s}^{\iy} p[x] v_1[x] \, dx\right)\right) (\cosh[\mu(s)]-2 \sinh[\mu(s)])+\\
\end{equation*}
\begin{equation*}
\left.c^2 (\int_{s}^{\iy} q[x] u[x] \, dx)^2 (2 \cosh[\mu(s)]-\sinh[\mu(s)])\right)+\\
\end{equation*}
\begin{equation*}
\left.\left.20 c \left(-\int_{s}^{\iy} q[x] v[x] \, dx+\int_{s}^{\iy} q_1[x] \, dx\right) \sinh[\mu(s)]\right)\right) n^{-\frac{2}{3}} +O(n^{-1})
\end{equation*}
\end{thm}

\clearpage \vspace{3ex} \noindent\textbf{\large Acknowledgements: }
The author would like to thank Alice and Aimee Choup for their continued affection and the Department of Mathematical Sciences at the
University of Alabama in Huntsville.

\end{document}